\documentclass[12pt]{article}
\usepackage{setspace} 
\usepackage[dvips]{graphicx} 
\usepackage{amsmath,amsthm,amsfonts}
\RequirePackage[colorlinks,citecolor=blue,urlcolor=blue,linkcolor=blue]{hyperref}
\hypersetup{
colorlinks = true,
citecolor=blue,
urlcolor=blue,
linkcolor=blue,
pdfauthor = {Alexey Kuznetsov},
pdfkeywords = {60G51,  exponential functional, Levy process, rational transform, Meijer G-function, Asian option},
pdftitle = {On the distribution of exponential functionals for L\'evy processes with jumps of rational transform},
pdfpagemode = UseNone
}
    \def\qed{\hfill$\sqcap\kern-8.0pt\hbox{$\sqcup$}$\\}
    \def\beq{\begin{eqnarray}}
    \def\eeq{\end{eqnarray}}
    \def\beqq{\begin{eqnarray*}}
    \def\eeqq{\end{eqnarray*}}
    
    \def\f{{\cal F}}
    \def\ee{\textnormal {e}}
    \def\re{\textnormal {Re}}
    \def\im{\textnormal {Im}}
    \def\p{{\mathbb P}}
    \def\mm{{\mathcal M}}
    \def\e{{\mathbb E}}
    \def\r{{\mathbb R}}
    \def\c{{\mathbb C}}
    \def\gg{{\mathcal G}}
    	
    \def\d{{\textnormal d}}
    \def\i{{\textnormal i}}

	\newtheorem{theorem}{Theorem}
	
	\newtheorem{proposition}{Proposition}
	\newtheorem{corollary}{Corollary}

\newtheorem{assumptionletter}{Assumption} 
 
\title{On the distribution of exponential functionals for L\'evy processes with jumps of rational transform}
\author{A. Kuznetsov
  \\ \\ 
Dept. of Mathematics and Statistics\\  York University
\\Toronto, ON, M3J 1P3 \\  Canada \\ e-mail: kuznetsov@mathstat.yorku.ca
 }\date{}

\begin{document}
\maketitle

\begin{abstract}
\bigskip
We derive explicit formulas for the Mellin transform and the distribution of the exponential functional for L\'evy processes with rational Laplace exponent. This extends recent results by Cai and Kou \cite{CaiKou2010} on the processes with hyper-exponential jumps. 
\end{abstract}

{\vskip 0.5cm}
{\vskip 0.5cm}


\section{Introduction}\label{intro}

Assume that $X$ is a L\'evy process and $\ee(q)$ is an independent random variable, which is exponentially distributed with parameter $q>0$. The exponential functional of $X$ is defined as
\beq\label{def_Iq}
I_q=\int\limits_0^{\ee(q)} e^{X_t} \d t.
\eeq
We can extend this definition to the case $q=0$, if we interpret $\ee(0)\equiv +\infty$ and assume that the process $X$ drifts to $-\infty$. Our main 
object of interest is the probability density function of $I_q$, defined as
\beqq
p(x)=\frac{\d }{\d x} \p(I_q\le x), \;\;\; x>0. 
\eeqq

Exponential functional of a L\'evy
process is a very interesting object, which has many applications in such areas as self-similar Markov processes, branching processes and Mathematical Finance. See \cite{BYS} for an overview of this topic. There exists an extensive literature covering the asymptotic behavior of $p(x)$ as $x\to +\infty$ 
(see \cite{CKP,CP,Maulik2006,RI,RI1} and the references therein) or as $x \to 0^+$ (see  \cite{PA, CR}). At the same time, the distribution of the exponential functional is known explicitly only in a few special cases: when $X$ is a standard Poisson process,  Brownian motion with drift, a particular spectrally negative Lamperti-stable process (see for instance \cite{CKP, KP,PP}), spectrally positive L\'evy process satisfying the Cram\'er's condition (see for instance \cite{PP1}). Note that in all above cases we have processes with one-sided jumps, and until very recent time 
there were no known examples of processes with double-sided jumps, for which the distribution of the exponential functional is known explicitly. 
However this situation has changed in the last several years.

First of all, in a recent paper \cite{CaiKou2010} (which is the main inspiration for our current work), Cai and Kou have obtained an explicit formula for the Mellin transform of 
$I_q$ in the case when $X$ has hyper-exponential jumps. While the authors did not provide a formula for the probability density function $p(x)$, 
it does follow from their results rather easily with the help of the theory of Meijer's G-function. Second, in \cite{KuzPardo2010} Kuznetsov and Pardo have derived 
an explicit infinite series representation and complete asymptotic expansions for $p(x)$ in the case of hypergeometric L\'evy processes. 

In this paper we will pursue the following three goals. First of all, we will generalize the results of Cai and Kou \cite{CaiKou2010} to the class of 
L\'evy processes with jumps of rational transform, and we will also cover the cases when $X$ has zero Gaussian component and/or zero linear drift. 
Second, we will show how to use the verification technique which was developed in \cite{KuzPardo2010}; we feel that
this technique is of independent interest as it considerably simplifies the derivation of many results on exponential functionals. 
Finally, we will present explicit formulas and complete asymptotic expansions for $p(x)$, which allow very simple numerical 
evaluation of this function and  other related quantities, such as the price of an Asian option with exponentially distributed maturity.

\section{Results}\label{results}

We will work with a L\'evy process $X$ with linear drift and diffusion coefficients $\mu\in \r$ and $\sigma \ge 0$ respectively, and for which the density of the L\'evy measure is given by 
\begin{equation}\label{def_pi_rational}
\pi(x)={\mathbf 1}_{\{x>0\}} \sum\limits_{j=1}^J 
\sum\limits_{i=1}^{m_j} \alpha_{ij} x^{i-1} e^{-\rho_j x}+{\mathbf 1}_{\{x<0\}}  \sum\limits_{j=1}^{\hat J} 
\sum\limits_{i=1}^{\hat m_j} \hat \alpha_{ij} x^{ i-1} e^{\hat \rho_j x}.
\end{equation}
We assume that $m_j, \hat m_j \in {\mathbb N}$, $\re(\rho_j)>0$ and $\re(\hat \rho_j)>0$ and that $\rho_i\ne \rho_j$ and  $\hat \rho_i \ne \hat \rho_j$ for all $i \ne j$. It is clear that $X$ has jumps of finite intensity $\lambda$, where
\beqq
\lambda=\int\limits_{\r} \pi(x) \d x=\sum\limits_{j=1}^J \sum\limits_{i=1}^{m_j} \alpha_{ij}  (i-1)!  \rho_j^{-i} 
+\sum\limits_{j=1}^{\hat J} \sum\limits_{i=1}^{\hat m_j} \hat \alpha_{ij} ( i-1)! \hat \rho_j^{-i}, 
\eeqq
in particular the process $X$ can be identified as
\beqq
X_t=\sigma W_t + \mu t +\sum\limits_{i=1}^{N_{t}} \xi_i,
\eeqq
where $N$ is a Poisson process with intensity $\lambda$ and $\xi_i$ are iid random variables, independent of $N$ and with distribution $\p(\xi_i \in \d x)=\lambda^{-1} \pi(x) \d x$.

Formula (\ref{def_pi_rational}) and the L\'evy-Khinchine formula imply that the Laplace exponent $\psi(z)=\ln \e [ \exp(z X_1)]$ is a rational function, which has the following partial fraction decomposition
\begin{equation}\label{rational_Laplace_exponent}
 \psi(z)=\frac{\sigma^2}2 z^2 + \mu z + \sum\limits_{j=1}^J \sum\limits_{i=1}^{m_j}  \frac{\alpha_{ij} (i-1)! }{(\rho_j-z)^{i}}
+\sum\limits_{j=1}^{\hat J} \sum\limits_{i=1}^{\hat m_j} \frac{\hat \alpha_{ij} (i-1)!}{  (\hat \rho_j+z)^{i}}  - \lambda.
\end{equation}
We see that $\psi(z)$ has poles at $\{\rho_j\}_{1\le j \le J}$ and $\{-\hat \rho_j\}_{1\le j \le \hat J}$ with the corresponding 
multiplicities $m_j$ and $\hat m_j$. Note that $\psi(z)$ has  $M=\sum_{j=1}^J m_j$ poles (counting with multiplicity) in the half-plane $\re(z)>0$ and
 $\hat M=\sum_{j=1}^{\hat J} \hat m_j$ poles in the half-plane $\re(z)<0$. 

The poles of $\psi(z)$ are important objects for our further results, but even more important are the zeros of $\psi(z)-q$, or, equivalently, the solutions
to the equation $\psi(z)=q$, where $q> 0$. First of all, let us find the total number of these solutions.  
Since $\psi(z)$ is a rational function of the form (\ref{rational_Laplace_exponent}), we can rewrite it as
\beq\label{eqn_psi_rational}
 \psi(z) = \frac{{\mathcal Q}(z)}{{\mathcal P}(z)}.
\eeq
The degree of the numerator \{denominator\} will be denoted as $Q={\textnormal {deg}}({\mathcal Q})$ \{$P={\textnormal {deg}}({\mathcal P})$\}. 
Note that (\ref{rational_Laplace_exponent}) implies  that $P=M+\hat M$ and 
 \beqq
\begin{cases}
  & Q=P+2  {\textnormal{ if }} \sigma>0, \\
  & Q=P+1 {\textnormal{ if }} \sigma=0  {\textnormal{ and }}  \mu\ne 0, \\
  & Q=P    {\textnormal{ if }} \sigma=\mu=0. 
\end{cases}
\eeqq

We claim that in all three cases the equation $\psi(z)=q$ has exactly $Q$ solutions $\zeta=\zeta(q)$ in $\c$ (counting multiplicity).  
This is obvious in the case 
$\sigma>0$ or $\sigma=0$ and $\mu\ne 0$, as in this case the equation $\psi(z)=q$ is equivalent to the polynomial equation ${\mathcal Q}(z)-q {\mathcal P}(z)=0$, which 
has degree $Q$. This statement is also true in the case $\sigma=\mu=0$: here we have $\psi(z)={\mathcal Q}(z)/{\mathcal P}(z) \to -\lambda$ as $z\to \infty$, therefore the leading term of the polynomial  ${\mathcal Q}(z)-q {\mathcal P}(z)$ can not be zero, thus the degree of this polynomial is also equal to $Q$ (note that the degree of this
polynomial is equal to $Q-1$ if $q=-\lambda<0$).

Let us denote the zeros of $\psi(z)-q$ in the half-plane $\re(z)>0$ \{$\re(z)<0$\} as $\zeta_1,\zeta_2,\dots,\zeta_K$   
\{$-\hat \zeta_1, -\hat \zeta_2,\dots, -\hat \zeta_{\hat K}$\}. Note that by definition we have $Q=K+\hat K$. We also assume that 
these numbers are labelled in the order of increase
of the real part, so that $\re(\zeta_j)\le \re(\zeta_{j+1})$ and 
$\re(\hat \zeta_j) \le  \re(\hat \zeta_{j+1})$. 

Let us summarize the definition of some important quantities, which will be used frequently later:
\begin{itemize}
 \item $M$ \{$\hat M$\} is the number of poles of $\psi(z)$ in the half-plane $\re(z)>0$ \{$\re(z)<0$\},   
 \item $K$ \{$\hat K$\} is the number of zeros of $\psi(z)-q$ in the half-plane $\re(z)>0$ \{$\re(z)<0$\},
 \item $P=M+\hat M$ \{$Q=K+\hat K$\} is the total number of poles \{zeros\} of $\psi(z)-q$.
\end{itemize}

We collect all important information about the solutions of the equation $\psi(z)=q$ in the next Proposition. 
\begin{proposition}\label{prop_roots_rational} Assume that $q>0$. Then
\indent
\begin{itemize}
 \item[(i)] $\;$  $\zeta_1$ and $\hat \zeta_1$ are real positive numbers. Moreover, for all $j\ne 1$, $\zeta_1<\re(\zeta_j)$ and $\hat \zeta_1< \re(\hat \zeta_j)$.
 \item[(ii)] $\;$  If $\sigma>0$ then $K=M+1$ and $\hat K=\hat M+1$. 
 \item[(iii)] $\;$ If $\sigma=0$, and $\mu>0$ \{$\mu<0$\} then $K= M+1$ and $\hat K=\hat M$ \{ $K=M$ and $\hat K=\hat M+1$\}.
 \item[(iv)] $\;$ If $\sigma=\mu=0$, then $K=M$ and $\hat K=\hat M$. 
 \item[(v)] $\;$ There exist at most $2Q-1$ complex numbers $q$ such that equation $\psi(z)=q$ has solutions of multiplicity greater than one.
 \item[(vi)] $\;$ As $q\to 0^+$ we have the following possibilities:
 \beqq
\begin{cases}
  & {\textnormal{if }} \e[X_1]>0, {\textnormal{ then }}   \zeta_1(0^+)=0 {\textnormal{ and }} \hat \zeta_1(0^+)>0, \\
  & {\textnormal{if }} \e[X_1]<0, {\textnormal{ then }}   \zeta_1(0^+)>0 {\textnormal{ and }} \hat \zeta_1(0^+)=0, \\
  & {\textnormal{if }} \e[X_1]=0, {\textnormal{ then }}   \zeta_1(0^+)=0 {\textnormal{ and }} \hat \zeta_1(0^+)=0.
\end{cases}
 \eeqq
\end{itemize}
\end{proposition}
\begin{proof}
 The proof of (i)-(iv) follows from Lemma 1.1 in \cite{Mordecki} 
(a somewhat simpler proof of a similar result has appeared recently in \cite{Fourati2010}).
The proof of  (vi) follows easily from the fact that $\psi(z)\sim \e[X_1] z + O(z^2)$ as $z\to 0$.  

Let us prove (v). Assume that equation $\psi(z)=q$ has a solution $z=z_0$
of multiplicity greater than one, therefore $\psi'(z_0)=0$. Using (\ref{eqn_psi_rational}) we find that $\psi'(z_0)=0$ implies ${\mathcal Q}'(z_0){\mathcal P}(z_0)-{\mathcal Q}(z_0){\mathcal P}'(z_0)=0$. 
The polynomial $H(z)={\mathcal Q}'(z){\mathcal P}(z)-{\mathcal Q}(z){\mathcal P}'(z)$ has degree at most $Q+P-1\le 2Q-1$, thus there exist at most $2Q-1$ distinct points $z_k$ for which $\psi'(z_k)=0$, 
which implies that there exist at most $2Q-1$ values $q$, given by $q_k=\psi(z_k)$, for which the equation $\psi(z)=q$ has solutions of multiplicity greater than one.
\end{proof}

{\bf Remark 1.}
Proposition \ref{prop_roots_rational}(v) states that in general it is very unlikely for the equation $\psi(z)=q$ to have multiple solutions, unless $q=0$ and 
$\e[X_1]=0$, in which case we have a double root at zero. In particular, when doing numerical computations, we can safely assume that
all solutions to $\psi(z)=q$ are simple. 
\\

 Our first goal is to identify the Mellin transform of $I_q$, that is $\e[I_q^{s-1}]$. Our main tool is the following result, which we have borrowed from \cite{KuzPardo2010}. Note that this result is not restricted to processes with jumps of rational transforms. 

\begin{proposition}\label{Ms_uniqueness_lemma} {\bf (Verification result)}
Assume that Cram\'er's condition is satisfied: there exists $z_0>0$ such that the Laplace exponent $\psi(z)$ is finite for all $z\in (0, z_0)$ and $\psi(\theta)=q$ for some $\theta \in (0,z_0)$. If $f(s)$ satisfies the following three properties
\begin{itemize}
 \item[(i)] $f(s)$ is analytic and zero-free in the strip $\re(s)\in (0,1+\theta)$,
 \item[(ii)] $f(1)=1$ and $f(s+1)=s f(s)/(q-\psi(s))$ for all $s\in (0,\theta)$,  
 \item[(iii)] $|f(s)|^{-1}=o(\exp(2 \pi |\im(s)|))$ as $\im(s)\to \infty$, $\re(s)\in (0,1+\theta)$,
\end{itemize}
 then $\e[I_q^{s-1}] \equiv f(s)$ for  $\re(s) \in (0,1+\theta)$.
\end{proposition}
\begin{proof}
The proof can be found in \cite{KuzPardo2010}, however for sake of completeness we will present the main steps of the proof here. 
Let us denote the Mellin transform of $I_q$ as $ \mm_Y(s)=\e[I_q^{s-1}]$. First of all, the Cram\'er's condition and 
Lemma 2 in \cite{r2007} imply that $\mm_Y(s)$ can be extended to an analytic function in the vertical strip $\re(s) \in (0,1+\theta)$. Since 
$|\mm_Y(s)|<\mm_Y(\re(s))$ we see that $\mm_Y(s)$ is bounded in the strip $\re(s) \in [\theta/2,1+\theta/2]$. It is well-known
(see  Lemma 2.1 in \cite{Maulik2006} or Proposition 3.1 in \cite{CarPetYor1997}) that $\mm_Y(s)$ satisfies the same functional equation as $f(s)$. Therefore the ratio $F(s)=\mm_Y(s)/f(s)$ is a periodic function: $F(s+1)=F(s)$; and due to condition (i) $F(s)$ can be extended to an analytic function
 in the entire complex plane. Finally, condition
(iii) and boundedness of $\mm_Y(s)$ imply that $F(s)=o(\exp(2\pi |\im(s)|)$ in the entire complex plane, and any function which is analytic, periodic with period equal to one, 
and which satisfies this upper bound must be identically equal to a constant. Since $F(1)=1$ we conclude that $F(s)\equiv 1$, that is $\mm_Y(s)\equiv f(s)$. 
\end{proof}


Proposition \ref{Ms_uniqueness_lemma} is a convenient tool which allows us to explicitly identify the Mellin transform of $I_q$ as a solution to the functional equation $f(s+1)=s f(s)/(q-\psi(s))$. What makes this equation analytically tractable is the fact that in our case
$s/(q-\psi(s))$ is a rational function. In fact, as was pointed out by the referee, functional equations of the form $f(s+1)=R(s)f(s)$ with a rational function $R(s)$ were studied by Mellin in a 1910 paper \cite{Mellin}, where he has also
investigated their connections with Gamma functions, Mellin transform and generalized hypergeometric functions.

\begin{theorem}\label{thm_main} Assume that $q>0$ or $q=0$ and $\e[X_1]<0$. For $\re(s) \in (0,1+\zeta_1)$ we have 
\beq\label{eq_main}
\e \left[ I_q^{s-1} \right]=A^{1-s} \times \Gamma(s) \times \frac{{\mathcal G}(s)}{{\mathcal G}(1)},
\eeq
where 
\beq\label{def_Bus}
{\mathcal G}(s)=\frac{\prod\limits_{j=1}^K \Gamma(1+\zeta_j-s)}{\prod\limits_{j=1}^J \Gamma(1+\rho_j-s)^{m_j} }
\times
\frac{\prod\limits_{j=1}^{\hat J} \Gamma(\hat \rho_j+s)^{\hat m_j}}{\prod\limits_{j=1}^{\hat K} \Gamma(\hat\zeta_j+s) },
\eeq
and the constant $A$ is defined as 
\beq\label{def_A}
\begin{cases}
   & A=\frac{\sigma^2}2, {\textnormal{ if }} \sigma>0, \\
   & A=|\mu|, {\textnormal{ if }} \sigma=0  {\textnormal{ and }} \mu \ne 0, \\
   & A=q+\lambda, {\textnormal{ if }} \sigma=\mu=0. \\
\end{cases}
\eeq
\end{theorem}
\begin{proof}
 Let us prove this Theorem in the case $\sigma=\mu=0$, the proof in other cases is similar. First of all, we check that the Cram\'er's condition
 in Proposition \ref{Ms_uniqueness_lemma} is satisfied
with $\theta=\zeta_1$. Let $f(s)$ be the function in the right-hand side of (\ref{eq_main}). 
Since $\Gamma(s)$ has no zeros in $\c$ and simple poles at $s \in \{0,-1,-2,\dots\}$, it is clear that $f(s)$ is analytic and zero-free in the vertical strip $\re(s) \in (0,1+\zeta_1)$. 

Next, we factorize the rational function $s/(q-\psi(s))$ as follows
\beq\label{s_q_psi_factorization}
\frac{s}{q-\psi(s)}=\frac{s}{q+\lambda}\frac{\prod\limits_{j=1}^J (s-\rho_j)^{m_j} }{\prod\limits_{j=1}^K (s-\zeta_j)}
\times
\frac{\prod\limits_{j=1}^{\hat J} (s+\hat \rho_j)^{\hat m_j}}{\prod\limits_{j=1}^{\hat K} (s+\hat\zeta_j) }.
\eeq
One can check that this is the correct factorization, since the rational functions in the right/left hand sides of this equation have identical
roots/poles and the same asymptotic behavior as $s\to \infty$ (recall that $\psi(s)\to -\lambda$ when $s\to \infty$ in the case $\sigma=\mu=0$).
Using (\ref{s_q_psi_factorization}) and the functional identity $\Gamma(s+1)=s\Gamma(s)$ we check that $f(s)$ satisfies the functional equation $f(s+1)=sf(s)/(q-\psi(s))$ for $\re(s) \in (0,1+\zeta_1)$.
Finally, we use the asymptotic expansion
\beq\label{gamma_asymptotics}
|\Gamma(x+\i y)|=\exp\left(-\frac{\pi}2|y|+\left(x-\frac12\right) \ln|y|+O(1) \right), \;\;\; y\to \infty 
\eeq
(see formula 8.328.1 in \cite{Jeffrey2007}), the fact that
\beqq
M=\sum\limits_{j=1}^J m_j=K, \;\;\; \hat M=\sum\limits_{j=1}^{\hat J} \hat m_j=\hat K
 \eeqq
(see Proposition \ref{prop_roots_rational}(iv)) and conclude
\beqq
|f(s)|=\exp\left(-\frac{\pi}2 |\im(s)|+O(\ln|s|)\right), \;\;\; \im(s) \to \infty
\eeqq
in the vertical strip $\re(s) \in (0,1+\zeta_1)$. Thus we see that $f(s)$ satisfies all three conditions of Proposition 
\ref{Ms_uniqueness_lemma}, therefore (\ref{eq_main}) is true for $\re(s) \in (0,1+\zeta_1)$. 
\end{proof}

{\bf Remark 2.}
Using Proposition \ref{prop_roots_rational} we see that in all cases, except when $\sigma=0$ and $\mu<0$, the right-hand side of (\ref{eq_main})
has more gamma functions in the numerator than in the denominator. This fact and the asymptotic formula (\ref{gamma_asymptotics})  imply that
$\e[I_q^{s-1}]$ decreases to zero exponentially fast as $\im(s) \to \infty$. In the case $\sigma=0$ and $\mu<0$ we have the same number of gamma functions 
in the numerator and denominator, therefore (with the help of (\ref{gamma_asymptotics})) we conclude that the rate of decrease is polynomial. This implies (via the inverse Mellin transform) that $p(x)$ is a smooth function on $\r^+$, except in the case  when $\sigma=0$ and $\mu<0$. This result is consistent with Proposition 2.1 in \cite{CarPetYor1997} (see also Proposition \ref{prop_p(x)} below).
\\

Formula (\ref{eq_main}) gives us the Mellin transform of $I_q$, which uniquely characterizes the distribution of $I_q$. It turns out that with little 
extra work we can obtain an explicit formula for the joint transform of $(X_{\ee(q)}, I_q)$. 
\begin{corollary}
 Assume that $q>0$. For  $\re(u) \in (-\hat \zeta_1,\zeta_1)$ and  $\re(s) \in (0,1+\zeta_1-\re(u))$ we have
\beq\label{eqn_2d_transform}
\e \left[ e^{uX_{\ee(q)}} I_q^{s-1} \right]=\frac{q A^{1-s}}{q-\psi(u)}  \times \Gamma(s) \times  \frac{{\mathcal G}(s+u)}{{\mathcal G}(1+u)}.
\eeq
\end{corollary}
\begin{proof}
 Assume that $u \in (-\hat \zeta_1, \zeta_1)$ and define a new measure $\tilde \p$ as the Escher transform of $\p$
\beqq
 \frac{\d \tilde \p}{\d \p}\bigg |_{\f_t}=e^{uX_t-t \psi(u)}. 
\eeqq
Define $I(t)=\int_0^t \exp(X_v)\d v$. 
Then we have
\beq\label{cor1_proof1}
\e\left[ e^{uX_{\ee(q)}} I_q^{s-1} \right]&=&q \int\limits_0^{\infty} e^{-qt} \e\left[ e^{uX_{t}} I(t)^{s-1} \right] \d t=
q \int\limits_0^{\infty} e^{-qt+t\psi(u)} \tilde\e\left[ I(t)^{s-1} \right] \d t=\frac{q}{\tilde q} \tilde\e\left[ I_{\tilde q}^{s-1} \right],
\eeq
where we have denoted $\tilde q=q-\psi(u)$. Note that $\psi(u)<0$ for $u\in (-\hat \zeta_1, \zeta_1)$, therefore $\tilde q>0$. The Laplace exponent of $X$ under the measure $\tilde \p$ is equal to $\tilde \psi(z)=\psi(z+u)-\psi(u)$,
therefore the zeros and poles of the rational function $\tilde\psi(z)-\tilde q$ can be obtained from the corresponding zeros and poles of $\psi(z)-q$ 
by shifting by $u$ units:
\beqq
\zeta_j \mapsto \zeta_j-u, \;\;\; \rho_j \mapsto \rho_j-u, \;\;\;
 {\hat \zeta_j} \mapsto \hat\zeta_j+u, \;\;\;  \hat \rho_j \mapsto \hat \rho_j+u, \;\;\;
\eeqq
We apply results of Theorem \ref{thm_main} and conclude that
\beq\label{cor1_proof2}
\tilde\e\left[ I_{\tilde q}^{s-1} \right]=A^{1-s}  \times \Gamma(s) \times  \frac{{\mathcal G}(s+u)}{{\mathcal G}(1+u)}.
\eeq
Identity (\ref{eqn_2d_transform}) for real values of $u$ follows from (\ref{cor1_proof1}) and (\ref{cor1_proof2}), the general case follows by analytic continuation.
\end{proof}

As we have mentioned above, formula (\ref{eq_main}) uniquely characterizes the distribution of $I_q$ via the inverse Mellin transform
\beq\label{p_inverse_Mellin}
p(x)=\frac{1}{2\pi \i} \int\limits_{1+\i \r} \e \left[ I_q^{s-1} \right] x^{-s} \d s. 
\eeq
In fact, $p(x)$ can be computed explicitly, and this can be achieved in a variety of ways. One approach (which is quite general) is to 
use the fact that the right-hand side of (\ref{eq_main}) is a meromorphic function in $\c$, then compute its residues, and then by shifting 
the contour of integration we will obtain either convergent series representations or complete asymptotic expansions for $p(x)$. This is the method that was 
used in \cite{KuzPardo2010} in the case of hypergeometric L\'evy processes. The second approach, which we will follow in this paper, 
is to use the theory of the Meijer G-function.

First of all, we want to ensure that the Mellin transform of $I_q$ given in (\ref{eq_main}) does not have multiple poles. One can still obtain explicit formula in 
the case when we do have multiple poles, however in this case we will have logarithmic terms and more lengthy expressions. We would like to 
avoid this situation, therefore from now on we will work under the following assumption, which ensures that $\Gamma(s) {\mathcal G(s)}$ does not have multiple poles.
\begin{assumptionletter} \label{as_A} 
         \renewcommand{\theenumi}{A.\arabic{enumi}} 
         \renewcommand{\labelenumi}{\theenumi} 
	${}$ \\
         \begin{enumerate} 
                 \item \label{as_A1} $\psi(z)-q$ has no multiple poles in the half-plane $\re(z)<0$.  
                 \item \label{as_A2} For $1\le i \le \hat M$ we have $\hat \rho_i \notin {\mathbb N}$.  
                 \item \label{as_A3} For $1\le i < j \le \hat M$ we have $\hat \rho_j-\hat \rho_i \notin {\mathbb N}$.  
                 \item \label{as_A4} $\psi(z)-q$ has no multiple zeros in the half-plane $\re(z)>0$.
                 \item \label{as_A5} For $1\le i < j \le K$ we have $\zeta_j-\zeta_i \notin {\mathbb N}$.
         \end{enumerate} 
\end{assumptionletter} 

\vspace{0.2cm}
Note that Assumptions \ref{as_A1}, \ref{as_A2} and  \ref{as_A3} can be easily verified, since usually we know explicitly the L\'evy measure, and therefore
the values of $\hat \rho_j$. Assumption \ref{as_A1} is equivalent to requiring $\hat m_j=1$ for $j=1,2,\dots,\hat J$. 
Assumption \ref{as_A4} can be verified numerically, however due to Proposition \ref{prop_roots_rational}(v) we would expect 
that this assumption will be satisfied in numerical experiments, unless we have a specifically engineered counter-example. 
Similarly, condition \ref{as_A5} should be satisfied in numerical experiments, since a ``generic'' polynomial equation normally does not have 
solutions which differ by an integer number. Thus Assumptions \ref{as_A4} and \ref{as_A5} are not very restrictive.

Next, let us introduce some notations. For any $w\in \c$ and  $n \in {\mathbb N}$ we define the vector $[w]_n \in \c^n$  as
\beqq
[w]_n=[w,w,\dots,w].
\eeqq
We define vectors ${\mathbf a} \in \c^{P+1}$ and ${\mathbf b} \in \c^{Q}$  as
\beqq
{\mathbf a}&=&[1,[1-\hat \rho_1]_{\hat m_1}, [1-\hat \rho_2]_{\hat m_2}, \dots, [1-\hat \rho_{\hat J}]_{\hat m_{\hat J}},
[1+\rho_1]_{m_1}, [1+\rho_2]_{m_2}, \dots , [1+\rho_J]_{m_J}],\\
{\mathbf b}&=&[1+\zeta_1,1+\zeta_2,\dots,1+\zeta_K,1-\hat \zeta_1, 1-\hat \zeta_2,\dots,1-\hat \zeta_{\hat K} ].
\eeqq
For any ${\mathbf x} \in \c^n$ and $j \in \{1,2,\dots,n\}$ we will denote as $(\mathbf x)^{j} \in \c^{n-1}$ the vector obtained from ${\mathbf x}$  by deleting the $j$-th coordinate, that is
\beqq
(\mathbf x)^{j}=[x_1,x_2,\dots,x_{j-1},x_{j+1},\dots,x_n].
\eeqq   
Finally, for ${\mathbf x}\in \c^n$ and $w \in \c$, we define the vector $w+{\mathbf x} \in \c^n$ as follows  
\beqq
w+{\mathbf x} = [w+x_1,w+x_2,\dots,w+x_n].
\eeqq  

For ${\mathbf u} \in \c^p$ and ${\mathbf v}\in \c^q$ the generalized hypergeometric function is defined as a formal power series
\beqq
{_pF_q}\left[
\begin{array} {c}
 {\mathbf u} \\  {\mathbf v} \\ \end{array} 
\bigg | z \right]=\sum\limits_{n\ge 0} \frac{(u_1)_n\dots (u_p)_n}{(v_1)_n\dots (v_q)_n} \frac{z^n}{n!}.
\eeqq
where $(w)_n=w(w+1)\dots(w+n-1)$ is the Pochhammer symbol. Using the ratio test, it is easy to see that the domain of convergence of this 
series is  (a) $\c$ when $p\le q$ (b) $\{|z|<1\}$ when $p=q+1$ (c) $\{0\}$ when $p\ge q+2$. 
In the latter case we will interpret ${_pF_q}$ as an asymptotic expansion as $z\to 0$.

Finally, let us define the following two functions (or formal power series in the case if they do not converge):

\beq\label{def_F1}
&&G_1(x)=\sum\limits_{j=1}^K  \;\frac{\prod\limits_{\substack {1\le i \le K \\ i \ne j}} \Gamma(b_i-b_j) \prod\limits_{i=1}^{\hat M+1} \Gamma(1+b_j-a_i) }
 {\prod\limits_{i=K+1}^Q \Gamma(1+b_j-b_i) \prod\limits_{i=\hat M+2}^{P+1} \Gamma(a_i-b_j)}  \\ \nonumber
 &&\qquad\qquad\qquad\times x^{-b_j} \; {_{P+1}F_{Q-1}}   
\left[
\begin{array} {c}
1+b_j-{\mathbf a} \\ 1+b_j-({\mathbf b})^{j} \\ \end{array} 
\bigg | (-1)^{M-K} x^{-1} \right] 
\eeq
\beq
\label{def_F2}
&& G_2(x)=\sum\limits_{j=1}^{\hat M+1} \; \frac{\prod\limits_{\substack {1\le i \le \hat M+1 \\ i \ne j}} \Gamma(a_j-a_i) \prod\limits_{i=1}^{K} \Gamma(1+b_i-a_j) }
 {\prod\limits_{i=\hat M+2}^{P+1} \Gamma(1+a_i-a_j) \prod\limits_{i=K+1}^{Q} \Gamma(a_j-b_i)} 
\\ \nonumber  &&\qquad\qquad\qquad\times
 x^{1-a_j} \; {_{Q}F_{P}}   
\left[
\begin{array} {c}
1-a_j+{\mathbf b} \\ 1-a_j+({\mathbf a})^j \\ \end{array} 
\bigg | (-1)^{\hat M-\hat K-1} x \right]
\eeq

Our final result is the following proposition, which completely describes the distribution of the exponential functional. For the definition of  the Meijer's G-function see formula 9.301 in \cite{Jeffrey2007}.
\begin{proposition}\label{prop_p(x)} Assume that $q>0$ or $q=0$ and $\e[X_1]<0$. Then $p(x)$ can be expressed in terms of 
Meijer's G-function as follows
 \beq\label{p_as_G}
 p(x)=\frac{A}{\gg(1)} \times G_{P+1,Q}^{K,\hat M+1} \left[ (Ax)^{-1} \bigg | 
\begin{array} {c}
{\mathbf a} \\ {\mathbf b} \\ \end{array}  \right].
\eeq
Moreover, we have the following series representations and asymptotic expansions for $p(x)$:
\\
\noindent {\bf (i)} if $\sigma>0$ then 
\beqq
p(x)&=&\frac{\sigma^2}{2{\mathcal G}(1)} G_1\left(\frac{\sigma^2 x}{2} \right), \;\;\; x>0, \\
p(x)&\sim& \frac{\sigma^2}{2{\mathcal G}(1)} G_2\left(\frac{\sigma^2 x}{2} \right), \;\;\; x\to 0^+;
 \eeqq
\\
\noindent {\bf (ii)} if $\sigma=0$ and $\mu \ne 0$ then 
\beqq
p(x)&=&\frac{|\mu|}{{\mathcal G}(1)} G_1\left(|\mu| x \right), \;\;\;   x> |\mu|^{-1}, \\
p(x)&=& \frac{|\mu|}{{\mathcal G}(1)} G_2\left( |\mu| x  \right), \;\;\; x< |\mu|^{-1};
 \eeqq
\noindent {\bf (i)} if $\sigma=\mu=0$ then 
\beqq
p(x)&\sim&\frac{q+\lambda}{{\mathcal G}(1)} G_1\left( (q+\lambda) x \right), \;\;\; x \to +\infty , \\
p(x)&=& \frac{q+\lambda}{{\mathcal G}(1)} G_2\left( (q+\lambda) x \right), \;\;\; x>0.
 \eeqq
\end{proposition}
\begin{proof}
Formula \eqref{p_as_G} follows from the expression for $p(x)$ as the inverse Mellin transform (\ref{p_inverse_Mellin}), formula (\ref{eq_main}) and  the definition of the Meijer's G-function 
(formula 9.301 in \cite{Jeffrey2007}). All the convergent series representations presented in Proposition \ref{prop_p(x)} follow immediately from 9.303 and 9.304 in \cite{Jeffrey2007}, while the asymptotic expansions follow from Theorem 1 in \cite{Fields1972}.
\end{proof}

\section{Acknowledgements}

The research was supported by the Natural Sciences and Engineering Research Council of Canada. We would like to thank Andreas Kyprianou, Juan Carlos Pardo and Mladen Savov for many stimulating and insightful discussions. We are also grateful to an anonymous referee for careful reading of the paper, for providing many valuable comments and for bringing the paper \cite{Mellin} to our attention.


\end{document}